\documentclass[11pt]{article}
\usepackage{amsfonts}
\usepackage{amsmath}
\usepackage{amssymb}
\usepackage{graphicx}
\usepackage[dvipsnames,usenames]{color}
\usepackage[pdfstartview=FitH ]{hyperref}
\usepackage[hyperpageref]{backref}
\begin{document}

\baselineskip 18pt
\def\o{\over}
\def\e{\varepsilon}
\title{\Large\bf On\ \ the\ \ Estimate\ \ for\ \ a\ \ Mean\ \ Value\\
Relative\ \ to\ \ $\bf {4\o p}={1\o n_1}+{1\o n_2}+{1\o n_3}$}
\author{Chaohua\ \  Jia}
\date{}
\maketitle {\small \noindent {\bf Abstract.} For the positive integer
$n$, let $f(n)$ denote the number of positive integer solutions
$(n_1,\,n_2,\,n_3)$ of the Diophantine equation
$$
{4\o n}={1\o n_1}+{1\o n_2}+{1\o n_3}.
$$
For the prime number $p$, $f(p)$ can be split into $f_1(p)+f_2(p),$
where $f_i(p)(i=1,\,2)$ counts those solutions with exactly $i$ of
denominators $n_1,\,n_2,\,n_3$ divisible by $p.$

Recently Terence Tao proved that
$$
\sum_{p< x}f_1(p)\ll x\exp({c\log x\o \log\log x})
$$
with other results. In this paper we shall improve it to
$$
\sum_{p< x}f_1(p)\ll x\log^5x\log\log^2x.
$$
}

\vskip.4in
\noindent{\bf 1. Introduction}

For the positive integer $n$, let $f(n)$ denote the number of positive
integer solutions $(n_1,\,n_2,\,n_3)$ of the Diophantine equation
$$
{4\o n}={1\o n_1}+{1\o n_2}+{1\o n_3}.
$$
Erd\" os and Straus conjectured that for all $n\geq 2,\, f(n)>0.$ It
is still an open problem now although there are some partial results.

In 1970, R. C. Vaughan[5] showed that the number of $n< x$ for which
$f(n)=0$ is at most $x\exp(-c\log^{2\o 3}x),$ where $x$ is sufficiently
large and $c$ is a positive constant.

Recently Terence Tao[4] studied the situation in which $n$ is the prime
number $p.$ He gave lower bound and upper bound for the mean value of
$f(p).$ Precisely, he split $f(p)$ into $f_1(p)+f_2(p),$ where $f_i(p)
(i=1,\,2)$ counts those solutions with exactly $i$ of denominators
$n_1,\,n_2,\,n_3$ divisible by $p.$ He proved that
\begin{align}
x\log^2x\ll\sum_{p< x}f_1(p)\ll x\exp({c\log x\o \log\log x})
\end{align}
and
\begin{align}
x\log^2x\ll\sum_{p< x}f_2(p)\ll x\log^2x\log\log x,
\end{align}
where $p$ denotes the prime number, $x$ is sufficiently large and $c$
is a positive constant.

For the progress and some explanation on the estimate in (2), one can see [2].
In this paper we shall improve the upper bound in (1).

{\bf Theorem.} Let $p$ denote the prime number. Then for sufficiently
large $x$, we have
$$
\sum_{p< x}f_1(p)\ll x\log^5x\log\log^2x.
$$

Throughout this paper, let $p$ denote the prime number, $c$ denote
the positive constant, $p(n)$ be the least prime factor of $n$, $P(n)$
be the largest prime factor of $n$, $d(n)$ be the divisor function,
$\varphi(n)$ be the Euler totient function, $\Omega(n)$ be the number
of prime factors of $n$ with multiplicity.

\vskip.3in
\noindent{\bf 2. Some preliminaries}

{\bf Lemma 1.} Let
$$
g(x)=a_nx^n+\cdots+a_1x+a_0
$$
be the polynomial in integer coefficients, $G(n)$ be the number of
solutions to the congruence equation
$$
g(x)\equiv 0\ ({\rm mod}\,n).
$$
Then $G(n)$ is a multiplicative function.

One can see page 34 of [1].

{\bf Lemma 2.} Let $g(x)$ be the polynomial in integer coefficients.
If
$$
g(x)\equiv 0,\quad g'(x)\equiv 0\ ({\rm mod}\,p)
$$
have no common solution, then the number of solutions to
$$
g(x)\equiv 0\ ({\rm mod}\,p^l)
$$
is equal to that to
$$
g(x)\equiv 0\ ({\rm mod}\,p).
$$

One can see page 36 of [1].

{\bf Lemma 3.} For the fixed integer $l$, let $G(n)$ be the number
of solutions to the congruence equation
$$
4lx^2+1\equiv 0\ ({\rm mod}\,n).
$$
Then
$$
G(n)\leq d(n).
$$

{\bf Proof.} By Lemma 1, we know
$$
G(p_1^{l_1}\cdots p_s^{l_s})=G(p_1^{l_1})\cdots G(p_s^{l_s}).
$$
Write
$$
g(x)=4lx^2+1.
$$
Then $g'(x)=8lx.$ It is obvious that
$$
g(x)\equiv 0,\quad g'(x)\equiv 0\ ({\rm mod}\,p)
$$
has no common solution. Thus Lemma 2 claims
$$
G(p^l)=G(p).
$$
It is easy to see that the congruence equation
$$
4lx^2+1\equiv 0\ ({\rm mod}\,p)
$$
has at most two solutions. Therefore
$$
G(p^l)=G(p)\leq 2\leq d(p^l).
$$
The conclusion of Lemma 3 follows.

{\bf Lemma 4.} For $x\geq 2,$ we have
$$
\sum_{n\leq x}{d^2(n)\o n}\ll\log^4 x.
$$

{\bf Proof.} Theorem 2 in [3] asserts that
$$
\sum_{n\leq x}d^2(n)\ll x\log^3x.
$$
Then
\begin{align*}
\sum_{n\leq x}{d^2(n)\o n}&\leq\sum_{i\leq\log_2 x}\sum_{2^i\leq n<
2^{i+1}}{d^2(n)\o n}\\
&\ll\sum_{i\leq\log_2 x}i^3\\
&\ll\log^4 x.
\end{align*}

{\bf Lemma 5.} Let
$$
\Psi(x,\,y)=\sum_{\substack{n\leq x\\ P(n)\leq y}}1.
$$
Then for $x\geq 10$, we have
$$
\Psi(x,\,\log x\log\log x)\ll\exp({3\log x\o (\log\log x)^{1\o 2}}).
$$

This is Lemma 1 in [3].

{\bf Lemma 6.} The estimate
$$
\sum_{\substack{Z^{1\o 2}\leq n\\ P(n)\leq Z^{1\o r}}}{d^2(n)\o
n}\ll\exp(\sum_{p\leq Z}{d^2(p)\o p}-{r\o 10}\log r)
$$
holds true for $1\leq r\leq {\log Z\o \log\log Z}$ uniformly.

This is a special case of Lemma 4 in [3].

\vskip.3in
\noindent{\bf 3. The proof of Theorem}

According to the discussion in the beginning of section 3 of [4], in
order to estimate
$$
\sum_{p< x}f_1(p),
$$
it is enough to estimate
\begin{align}
\sum_{\substack{a,\,l\\ al\leq x}}{x d(4la^2+1)\o \varphi(4al)
\log(1+{x\o al})}.
\end{align}
Using the bound
$$
\varphi(n)\gg{n\o \log\log n},
$$
we should estimate
$$
x\log\log x\sum_{\substack{a,\,l\\ al\leq x}}{d(4la^2+1)\o al
\log(1+{x\o al})}
$$
or
\begin{align*}
x&\log\log x\sum_{i\leq\log_2x}\sum_{j\leq\log_2x-i}{1\o 1+\log_2x-i-j}\cdot\\
&\quad\cdot{1\o 2^{i+j}}\sum_{2^i< a\leq 2^{i+1}}\sum_{2^j< l\leq
2^{j+1}}d(4la^2+1).
\end{align*}

Now we consider the estimate for the sum
\begin{align}
\sum_{V< l\leq 2V}\sum_{W< a\leq 2W}d(4la^2+1).
\end{align}
We shall use some ideas from [3].

Firstly assume that $V\leq W.$  Let
\begin{align}
Z=W^{1\o 20}.
\end{align}
Write $n$ uniquely as
$$
n=p_1^{s_1}\cdots p_j^{s_j}\cdot p_{j+1}^{s_{j+1}}\cdots p_r^{s_r},
\qquad p_1<\cdots< p_j< p_{j+1}<\cdots< p_r,
$$
where
$$
p_1^{s_1}\cdots p_j^{s_j}\leq Z< p_1^{s_1}\cdots p_j^{s_j}p_{j+1}^{s_{j+1}}.
$$
We can decompose $4la^2+1$ as
\begin{align}
4la^2+1=(p_1^{s_1}\cdots p_j^{s_j})(p_{j+1}^{s_{j+1}}\cdots p_r^{s_r})
=b(l,\,a)c(l,\,a),
\end{align}
where
$$
b(l,\,a)\leq Z,\qquad\quad (b(l,\,a),\,c(l,\,a))=1.
$$
We shall discuss in four cases as in [3].

Case I.\ $p(c(l,\,a))>Z^{1\o 2}.$

Since $p(c(l,\,a))>Z^{1\o 2},\,d(c(l,\,a))=O(1)$. Thus
$$
d(4la^2+1)=d(b(l,\,a))d(c(l,\,a))\ll d(b(l,\,a)).
$$
Lemmas 3 and 4 yield that
\begin{align*}
{\sum}_{{\rm I}}&\ll\sum_{b\leq Z}d(b)\sum_{V< l\leq 2V}\sum_{\substack{W<
a\leq 2W\\ 4la^2+1\equiv 0\,({\rm mod}\,b)}}1\\
&\ll\sum_{b\leq Z}d(b)\sum_{V< l\leq 2V}{W\o b}\sum_{\substack{a=1\\
4la^2+1\equiv 0\,({\rm mod}\,b)}}^b1\\
&\ll VW\sum_{b\leq Z}{d^2(b)\o b}\\
&\ll VW\log^4(2W).
\end{align*}

Case II.\ $p(c(l,\,a))\leq Z^{1\o 2},\,b(l,\,a)\leq Z^{1\o 2}.$

Write $p=p(c(l,\,a)).$ Then $p^s\| 4la^2+1,\,p\leq Z^{1\o 2}.$ The
fact that $b(l,\,a)\leq Z^{1\o 2},\,b(l,\,a)p^s>Z$ yields $p^s>
Z^{1\o 2}.$ Let $s_p$ be the smallest $s$ such that $p^s> Z^{1\o 2}.$
Thus $s_p\geq 2.$ On the other hand, $p^{s_p\o 2}\leq p^{s_p-1}\leq
Z^{1\o 2}\Longrightarrow p^{s_p}\leq Z.$ Now we have
$$
{1\o p^{s_p}}\leq\min({1\o Z^{1\o 2}},\,{1\o p^2}).
$$
Hence
\begin{align*}
\sum_{p\leq Z^{1\o 2}}{1\o p^{s_p}}&\leq\sum_{p\leq Z^{1\o 4}}{1\o
Z^{1\o 2}}+\sum_{Z^{1\o 4}< p}{1\o p^2}\\
&\ll Z^{-{1\o 4}}.
\end{align*}

Lemmas 3 yields that
\begin{align*}
{\sum}_{{\rm II}}&\ll W^\e\sum_{p\leq Z^{1\o 2}}\sum_{V< l\leq 2V}
\sum_{\substack{W< a\leq 2W\\ 4la^2+1\equiv 0\,({\rm mod}\,p^{s_p})}}1\\
&\ll W^\e\sum_{p\leq Z^{1\o 2}}\sum_{V< l\leq 2V}{W\o p^{s_p}}
\sum_{\substack{a=1\\ 4la^2+1\equiv 0\,({\rm mod}\,p^{s_p})}}^{p^{s_p}}1\\
&\ll VW^{1+\e}\sum_{p\leq Z^{1\o 2}}{1\o p^{s_p}}\\
&\ll VW^{1-{1\o 80}+\e}\ll VW,
\end{align*}
where $p^{s_p}\leq Z$ works.

Case III.\ $p(c(l,\,a))\leq \log W\log\log W,\,b(l,\,a)> Z^{1\o 2}.$

We have $p(c(l,\,a))\leq\log W\log\log W\Longrightarrow P(b(l,\,a))<
\log W\log\log W.$ Then Lemmas 3 and 5 yield that
\begin{align*}
{\sum}_{{\rm III}}&\ll W^\e\sum_{\substack{Z^{1\o 2}< b\leq Z\\ P(b)<
\log W\log\log W}}\sum_{V< l\leq 2V}\sum_{\substack{W< a\leq 2W\\
4la^2+1\equiv 0\,({\rm mod}\,b)}}1\\
&\ll W^\e\sum_{\substack{Z^{1\o 2}< b\leq Z\\ P(b)< \log W\log\log W}}
{d(b)\o b}\,VW\\
&\ll VW^{1+2\e}Z^{-{1\o 2}}\sum_{\substack{b\leq Z\\ P(b)< \log W\log\log W}}1\\
&\ll VW^{1-{1\o 40}+2\e}\Psi(W,\,\log W\log\log W)\\
&\ll VW^{1-{1\o 40}+3\e}\ll VW.
\end{align*}

Case IV.\ $\log W\log\log W< p(c(l,\,a))\leq Z^{1\o 2},\,b(l,\,a)>
Z^{1\o 2}.$

Let
$$
r_0=[{\log Z\o \log(\log W\log\log W)}].
$$
Since
$$
\log W\log\log W> Z^{1\o r_0+1},
$$
for $2\leq r\leq r_0$, we consider these $(l,\,a)$ which satisfy
$$
Z^{1\o r+1}< p(c(l,\,a))\leq Z^{1\o r}
$$
so that
$$
P(b(l,\,a))< p(c(l,\,a))\leq Z^{1\o r}.
$$

We have
$$
\Omega(c(l,\,a))\leq {3\log W\o \log p(c(l,\,a))}\leq {3(r+1)\log W\o \log Z}\\
\leq 60(r+1)\leq 120r
$$
so that
$$
d(c(l,\,a))\leq A^r,
$$
where $A$ is a positive constant.

Lemmas 3, 4 and 6 yield that
\begin{align*}
{\sum}_{{\rm IV}}&\ll\sum_{2\leq r\leq r_0} A^r\sum_{\substack{Z^{1\o
2}< b\leq Z\\ P(b)< Z^{1\o r}}}d(b)\sum_{V< l\leq 2V}\sum_{\substack{W<
a\leq 2W\\ 4la^2+1\equiv 0\,({\rm mod}\,b)}}1\\
&\ll VW\sum_{2\leq r\leq r_0} A^r\sum_{\substack{Z^{1\o 2}< b\leq Z
\\ P(b)< Z^{1\o r}}}{d^2(b)\o b}\\
&\ll VW\sum_{2\leq r\leq r_0} A^r\exp(\sum_{p\leq Z}{d^2(p)\o p}-{r\o 10}\log r)\\
&= VW\sum_{2\leq r\leq r_0} A^r\exp(\sum_{p\leq Z}{4\o p}-{r\o 10}\log r)\\
&\ll VW\sum_{2\leq r\leq r_0} A^r\exp(4\log\log Z-{r\o 10}\log r)\\
&\ll VW\log^4 Z\sum_{r=2}^\infty A^r\exp(-{r\o 10}\log r)\\
&\ll VW\log^4(2W),
\end{align*}
where the series is convergent.

Then assume that $W< V.$ In this situation, we shall change the role
of $l$ and $r$ and shall consider the linear congruence equation
$$
4a^2l+1\equiv 0\,({\rm mod}\,n)
$$
for the fixed $a$. This situation is simpler than previous one. We can
get the similar estimate as above.

Combining all of above, we get
\begin{align*}
\sum_{V< l\leq 2V}\sum_{W< a\leq 2W}d(4la^2+1)\ll VW\log^4(2W).
\end{align*}
Hence,
\begin{align*}
x&\log\log x\sum_{i\leq\log_2x}\sum_{j\leq\log_2x-i}{1\o 1+\log_2x-i-j}\cdot\\
&\quad\cdot{1\o 2^{i+j}}\sum_{2^i< a\leq 2^{i+1}}\sum_{2^j< l\leq
2^{j+1}}d(4la^2+1)\\
&\ll x\log^4 x\log\log x\sum_{i\leq\log_2x}\sum_{j\leq\log_2x-i}{1\o 1+\log_2x-i-j}\\
&\ll x\log^4 x\log\log x\sum_{i\leq\log_2x}\sum_{1\leq h\leq\log_2x-i+1}{1\o h}\\
&\ll x\log^4 x\log\log x\sum_{i\leq\log_2x}\log(\log_2x-i+2)\\
&\ll x\log^5 x\log\log^2 x.
\end{align*}

So far the proof of Theorem is finished.


\bigskip

\

Institute of  Mathematics, Academia Sinica, Beijing 100190, P. R.
China

E-mail: jiach@math.ac.cn
\end{document}